\documentclass[aop,MSNbibl,dvips]{arximspdf}
\usepackage{graphicx}
\doi{10.1214/12-AOP756} 
\volume{41}
\issue{5}
\pubyear{2013}
\firstpage{3241}
\lastpage{3260}

\makeatletter
\newcommand{\rrVert}{\Vert}
\newcommand{\rrvert}{\vert}
\newcommand{\llVert}{\Vert}
\newcommand{\llvert}{\vert}

\newcommand{\underset}[2]{\mathop{#2}_{#1}}

\newtheorem{theorem}{Theorem}[section]
\newtheorem{corollary}[theorem]{Corollary}
\newproclaim{definition}[theorem]{Definition}
\newtheorem{lemma}[theorem]{Lemma}
\newtheorem{proposition}[theorem]{Proposition}

\makeatother

\begin{document}
\begin{frontmatter}

\title{An inequality for the distance between densities of~free convolutions}
\runtitle{An inequality for free convolutions}

\begin{aug}
\author[A]{\fnms{V.} \snm{Kargin}\corref{}\ead[label=e1]{v.kargin@statslab.cam.ac.uk}}
\runauthor{V. Kargin}
\affiliation{University of Cambridge}
\address[A]{Statistical Laboratory\\
Department of Pure Mathematics\\
\quad and Mathematical Statistics\\
University of Cambridge, Cambridge\\
United Kingdom\\
\printead{e1}} 
\end{aug}

\received{\smonth{8} \syear{2011}}
\revised{\smonth{1} \syear{2012}}

%
\begin{abstract}
This paper contributes to the study of the free additive convolution of
probability measures. It shows that under some conditions, if measures
$\mu_{i}$ and $\nu_{i}$, $i=1,2$, are close to each other in terms of the
L\'{e}%
vy metric and if the free convolution $\mu_{1}\boxplus\mu_{2}$ is
sufficiently smooth, then $\nu_{1}\boxplus\nu_{2}$ is absolutely
continuous, and the densities of measures $\nu_{1}\boxplus\nu_{2}$
and $%
\mu_{1}\boxplus\mu_{2}$ are close to each other. In particular,
convergence in distribution $\mu_{1}^{ ( n ) }\rightarrow
\mu_{1}, $ $\mu_{2}^{ ( n ) }\rightarrow\mu_{2}$ implies that
the density of $\mu_{1}^{ ( n ) }\boxplus\mu_{2}^{ (
n ) }$ is defined for all sufficiently large $n$ and converges to
the density of $\mu_{1}\boxplus\mu_{2}$. Some applications are
provided, including: (i) a new proof of the local version of the free
central limit theorem, and (ii) new local limit theorems for sums of
free projections, for sums of $\boxplus$-stable random variables and
for eigenvalues of a sum of two $N$-by-$N$ random matrices.
\end{abstract}

%

\end{frontmatter}

\section{Introduction}

Free convolution is a binary operation on the set of probability
measures on
the real line that converts this set into a commutative semigroup. In
contrast to the usual convolution, this operation is nonlinear
relative to
taking convex combinations of measures. The study of properties of free
convolution is motivated by its numerous applications to operator
algebras~\cite{voiculescu86,voiculescu96,haagerupthorbjornsen05},
random matrices \cite
{voiculescu91,speicher93,pasturvasilchuk00,guionnetkrishnapurzeitouni11},
representations of the symmetric group~\cite{biane98} and quantum physics
\cite{zee96,feinbergzee98}.

Starting with work by Voiculescu~\cite{voiculescu86}, it was noted that
free convolution has strong smoothing properties. Let $\mu_{1}\boxplus
\mu_{2}$ denote the free convolution of probability measures $\mu_1$
and $\mu_2$. In~\cite{bercovicivoiculescu98}, it was proved that $\mu
_{1}\boxplus \mu_{2}$ has an atom at $x \in\mathbb{R}$ if and only if
there are $y\in \mathbb{R}$ and $z\in\mathbb{R}$ such that $x=y+z$, and
$\mu_{1} ( \{ y \} ) +\mu_{2} ( \{ z \} ) >1$. In~\cite{belinschi08},
it was shown that $\mu_{1}\boxplus\mu_{2}$ can have a singular
component if and only if one of the measures is concentrated on one
point, and the other has a singular component (so that the resulting
free convolution is simply a translation of the measure with the
singular component). Moreover, in the same paper it was shown that the
density of the absolutely continuous part of the free convolution
measure is analytic wherever the density is positive and finite.

Some quantitative versions of the smoothing property of free convolution
have also been given. In particular, in~\cite{voiculescu93} it was shown
that if $\mu_{1}$ is absolutely continuous with density $f_{\mu
_{1}}\in
L^{p} ( \mathbb{R} )$ ($p\in( 1,\infty]
$), then the
free convolution of $\mu_{1}$ with an arbitrary other measure $\mu_{2}$ is
absolutely continuous with density $f_{\mu_{1}\boxplus\mu_{2}}\in$ $
L^{p} ( \mathbb{R} )$, and $\llVert f_{\mu
_{1}\boxplus\mu
_{2}}\rrVert_{p}\leq\llVert f_{\mu_{1}}\rrVert_{p}$. In
particular, the supremum of the density $f_{\mu_{1}\boxplus\mu_{2}}$ is
less than or equal to the supremum of the density of $f_{\mu_{1}}$.

Another important property of free convolution is that it is continuous
with respect to weak convergence\vspace*{1pt} of measures. In
particular, by a result in~\cite{bercovicivoiculescu93}, if $\mu_{1}^{
( N ) }\rightarrow\mu_{1} $ and $\mu_{2}^{ ( N ) }\rightarrow\mu_{2}$
as $N$ grows to infinity (where $\rightarrow$ denotes convergence in
distribution),
then $%
\mu_{1}^{ ( N ) }\boxplus\mu_{2}^{ ( N )
}\rightarrow
\mu_{1}\boxplus\mu_{2}$. In fact, Theorem~4.13 in
\cite{bercovicivoiculescu93} says that $d_{L} ( \mu_{1}\boxplus\mu
_{2},\nu
_{1}\boxplus\nu_{2} ) \leq d_{L} ( \mu_{1},\nu_{1} ) +d_{L} ( \mu
_{2},\nu_{2} )$, where $d_{L}$ denotes
the L\'{e}vy distance on the set of probability measures on $\mathbb{R}$.

The main result of this paper establishes a strengthened version of this
property. If distances $d_{L} ( \mu_{1},\nu_{1} ) $ and $%
d_{L} ( \mu_{2},\nu_{2} ) $ are sufficiently small, and if
$\mu_{1}\boxplus\mu_{2}$ is sufficiently smooth, then $\nu_{1}\boxplus
\nu_{2}$ is absolutely continuous and the distance between the densities
of $%
\mu_{1}\boxplus\mu_{2}$ and $\nu_{1}\boxplus\nu_{2}$ can be
bounded in
terms of the L\'{e}vy distances between the original measures.

In particular, this result shows that the convergence in distribution
$\mu_{1}^{ ( N ) }\rightarrow\mu_{1}$ and $\mu_{2}^{ (
N )
}\rightarrow\mu_{2}$ implies the convergence of the probability densities
of $\mu_{1}^{ ( N ) }\boxplus\mu_{2}^{ ( N )
}$ to the
density of $\mu_{1}\boxplus\mu_{2}$.

We prove this result under an assumption imposed on the measures $\mu_{1}$
and $\mu_{2}$, which we call the smoothness of the pair $ ( \mu_{1},\mu
_{2} ) $ at a point of its support $x$. This assumption holds
at a generic point $x$ if $\mu_{1}=\mu_{2}=\mu$, and the density of
$\mu
\boxplus\mu$ is absolutely continuous and positive at $x$. In the case
when $\mu_{1}\neq\mu_{2}$, this assumption should be checked
directly. We
envision that in applications $\mu_{1}$ and $\mu_{2}$ are fixed measures
for which this assumption can be directly checked, and $\mu_{1}^{ (
N ) }$ and $\mu_{2}^{ ( N ) }$ are (perhaps random) measures
for which it can be checked that they are close to $\mu_{1}$ and $\mu_{2}$
in the L\'{e}vy distance.

In order to formulate our main result precisely, we introduce several
definitions. Let $\mu_{1}$ and $\mu_{2}$ be two probability measures
on $%
\mathbb{R}$ with the Stieltjes transforms $m_{\mu_{1}} ( z
) $ and
$m_{\mu_{2}} ( z )$, where the Stieltjes transform of a
probability measure $\mu$ is defined by the formula
\[
m_{\mu}(z):=\int_{\mathbb{R}}\frac{\mu(dx)}{x-z}.
\]
Then, the free convolution $\mu_{1}\boxplus\mu_{2}$ is defined as a
probability measure on $\mathbb{R}$ with the Stieltjes transform
$m_{\mu
_{1}\boxplus\mu_{2}} ( z ) $, which satisfies the following
system of equations:
%
%
\begin{eqnarray}\label{system}
m_{\mu_{1}\boxplus\mu_{2}} ( z ) &=&m_{\mu_{1}} \bigl( \omega_{1} ( z )
\bigr),
\nonumber\\
m_{\mu_{1}\boxplus\mu_{2}} ( z ) &=&m_{\mu_{2}} \bigl( \omega_{2} ( z )
\bigr),
\\
z-\frac{1}{m_{\mu_{1}\boxplus\mu_{2}} ( z ) } &=&\omega_{1} ( z )
+\omega_{2} ( z ).
\nonumber
\end{eqnarray}
Here $\omega_{1} ( z ) $ and $\omega_{2} ( z )
$ are
analytic functions in $\mathbb{C}^{+}:= \{ z\dvtx\Im z>0 \}$,
that map
$\mathbb{C}^{+}$ to itself, that have the property $\Im\omega_{j} (
z ) \geq\Im z$, and such that $\omega_{j}(z)=z+o(z)$ as
$z\rightarrow
\infty$ in the sector $\Im z>\kappa\llvert\Re z\rrvert$,
where $%
\kappa$ is an arbitrary positive constant~\cite{biane98b}.
Functions $%
\omega_{1} ( z ) $ and $\omega_{2} ( z ) $ are called
the \textit{subordination functions} for the pair $ ( \mu_{1},\mu_{2} )$.

The definition of free convolution by the system (\ref{system}) is
equivalent to the standard definition through $R$-transforms
(\cite{voiculescudykemanica92} and~\cite{nicaspeicher06}) if one sets
$\omega_{1} ( z ) =z-R_{\mu_{2}} ( -m_{\mu_{1}\boxplus\mu
_{2}} ( z ) ) $, and similarly for $\omega_{2} (
z )$.

The subordination functions $\omega_{j}(z)$ depend not only on $z$ but also
on the pair $ ( \mu_{1},\mu_{2} )$. In particular, some
properties of the measures $\mu_{1}$ and $\mu_{2}$ are encoded in the
functions $\omega_{j}$. A proper but more cumbersome notation would be
$%
\omega_{j}(\mu_{1},\mu_{2},z)$ where $j=1,2$. In the cases when we need
to compare the subordination functions for pairs $ ( \mu_{1},\mu_{2} )
$ and $ ( \nu_{1},\nu_{2} )$, we will denote
them by
$\omega_{\mu,j}(z)$ and $\omega_{\nu,j} ( z ) $, respectively.

The system (\ref{system}) implies the following system of equations
for $\omega_{j}$:
%
%
\begin{eqnarray}\label{systemfort}
\frac{1}{z-\omega_{1} ( z ) -\omega_{2} ( z ) } &=&m_{\mu_{1}} \bigl(
\omega_{1} ( z ) \bigr),
\nonumber\\[-8pt]\\[-8pt]
\frac{1}{z-\omega_{1} ( z ) -\omega_{2} ( z ) } &=&m_{\mu_{2}} \bigl(
\omega_{2} ( z ) \bigr).
\nonumber
\end{eqnarray}

Note that the analytic solutions of the system (\ref{systemfort}) that
satisfy the asymptotic condition at infinity are unique in $\mathbb{C}^{+}$.
(This follows from the facts that the solutions are unique in the area
$\Im
z\geq\eta_{0}$ for sufficiently large $\eta_{0}$ and that the analytic
continuation in a simply-connected domain is unique.)

By Theorem 3.3 in~\cite{belinschi08}, the limits $\omega_{j}(x)=\lim
_{\eta
\downarrow0}\Im\omega_{j} ( x+i\eta) $ exist, and we
make the
following definition.

%
\begin{definition}
A pair of probability measures on the real line $ ( \mu_{1},\mu_{2} ) $
is said to be \textit{smooth} at $x$ if the following two
conditions hold:\vspace*{8pt}

\mbox{\hphantom{i}}(i) $\Im\omega_{j}(x)>0$ for $j=1,2$, and

(ii)
%
%
\begin{equation}\label{genericity}
k_{\mu}(x):=\frac{1}{m_{\mu_{1}}^{\prime}(\omega_{1}(x))}+\frac
{1}{%
m_{\mu_{2}}^{\prime}(\omega_{2}(x))}-\bigl(x-
\omega_{1}(x)-\omega_{2}(x)\bigr)^{2}\neq0.
\end{equation}
\end{definition}

Inequality (\ref{genericity}) is a technical condition and holds for a
generic point $x\in\mathbb{R}$.\vadjust{\goodbreak}

Condition (i) is somewhat stronger than the condition that $\mu
_{1}\boxplus\mu_{2}$ is Lebesgue absolutely continuous at $x$.
Indeed, if $%
\Im\omega_{j}(x)>0$ for $j=1,2$, then the limit
\[
\lim_{\eta\rightarrow0}m_{\mu_{1}\boxplus\mu_{2}} ( x+i\eta) =\lim
_{\eta\rightarrow0}m_{\mu_{1}}
\bigl( \omega_{1} ( z ) \bigr)
\]
exists and is finite. By using results in~\cite{belinschi08}, we can infer
from this fact that $\mu_{1}\boxplus\mu_{2}$ is Lebesgue absolutely
continuous at $x$.

In the converse direction, we have only that if $\mu_{1}=\mu_{2}=\mu$,
and $\mu\boxplus\mu$ is absolutely continuous with positive density
at $%
x, $ then condition (i) in the definition of smoothness is satisfied;
see Proposition~\ref{propositionsmoothnessstability} below.

The fact that smoothness is strictly stronger than absolute continuity
of $%
\mu_{1}\boxplus\mu_{2}$ can be seen from the following example. If
$\mu_{1}$ is a point mass at $0$, that is, $\mu_{1}=\delta_{0}$, and if
$\mu_{2}$ is absolutely continuous at $x$, then $\mu_{1}\boxplus\mu
_{2}=\mu_{2}$ is absolutely continuous at $x$, but the pair $ ( \mu
_{1},\mu_{2} ) $ is not smooth at $x$. Indeed, $m_{\delta_{0}}=-z^{-1}$,
and system (\ref{systemfort}) implies that $\omega_{2} (
z ) =z$.
Hence, $\Im\omega_{2} ( x ) =0$ for every $x$.

On the other hand smoothness holds for many examples that we consider below.

Next, let us recall the following standard definition.

%
\begin{definition}
The \textit{L\'{e}vy distance} between probability measures $\mu$ and
$\nu$
is
\[
d_{L} ( \mu,\nu) =\sup_{x}\inf\bigl\{ s
\geq0\dvtx F_{\nu
} ( x-s ) -s\leq F_{\mu} ( x ) \leq F_{\nu}
( x+s ) +s \bigr\},
\]
where $F_{\mu} ( t ) $ and $F_{\nu} ( t ) $ are the
cumulative distribution functions of $\mu$ and $\nu$.
\end{definition}

It is well known that $\mu^{ ( N ) }\rightarrow\mu$ in
distribution (i.e., the cumulative distribution function of $\mu^{ (
N ) }$ weakly converges to the cumulative distribution function
of $\mu
$) if and only if $d_{L} ( \mu^{ ( N ) },\mu)
\rightarrow0;$ see, for example, Theorem III.1.2 on page 314 and Exercise
III.1.4 on page 316 in~\cite{shiryaev96}.

Here is the main result of this paper.

%
\begin{theorem}
\label{theoremmain}Assume that a pair of probability measures $(\mu
_{1},\mu_{2})$ is smooth at $x$. Then, there are some $s_{\mu,0}>$
and $%
c_{\mu}>0$, which depend only on $(\mu_{1},\mu_{2})$, such that for all
pairs of probability measures $ ( \nu_{1},\nu_{2} ) $ with
$%
d_{L} ( \mu_{j},\nu_{j} ) <s\leq s_{\mu,0}$ for both $j=1,2$,
it is true that $\nu_{1}\boxplus\nu_{2}$ is absolutely continuous in a
neighborhood of $x$, and
\[
\bigl|f_{\nu_{1}\boxplus\nu_{2}} ( x ) -f_{\mu_{1}\boxplus
\mu
_{2}}(x)\bigr|<c_{\mu}s,
\]
where $f_{\nu_{1}\boxplus\nu_{2}}$ and $f_{\mu_{1}\boxplus\mu_{2}}$
are the densities of $\nu_{1}\boxplus\nu_{2}$ and $\mu_{1}\boxplus
\mu_{2}$, respectively.
\end{theorem}

This theorem will be proved as a corollary to Proposition
\ref{propositiondifferenceoft} below. The assumptions of the theorem are
sufficient but possibly not necessary. Of course, it is necessary to require
that $\mu_{1}\boxplus\mu_{2}$ be absolutely\vadjust{\goodbreak} continuous at $x$ so
that the
density $f_{\mu_{1}\boxplus\mu_{2}} ( x ) $ is well
defined. In
addition, a simple example shows that absolute continuity alone is not
sufficient. Indeed, if $\mu_{1}=\nu_{1}=\delta_{0}$ is a point mass at
zero, and $\mu_{2}$ is absolutely continuous, then $\delta_{0}\boxplus
\mu_{2}$ is absolutely continuous, but $\delta_{0}\boxplus\nu_{2}$ is not
necessarily so, even if $\nu_{2}$ is close to $\mu_{2}$ in the L\'{e}vy
distance. However, it is not clear if the assumption of absolute continuity
of $\mu_{1}\boxplus\mu_{2}$ implies the statement of the theorem once
this degenerate case is ruled out.

The constant $c_{\mu}$ in the theorem can be bounded in terms of $\Im
\omega_{\mu,j} ( x ) $ and $\llvert k_{\mu}(x)\rrvert$
from (\ref{genericity}). In particular, if $\Im\omega_{\mu,j} (
x ) $ and $\llvert k_{\mu}(x)\rrvert$ are uniformly bounded
away from zero for all $x\in( a,b )$, then $\sup_{x\in
(
a,b ) }|f_{\nu_{1}\boxplus\nu_{2}} ( x ) -f_{\mu
_{1}\boxplus\mu_{2}}(x)|<cs$ for some $c>0$.

The main ideas of the proof of Theorem~\ref{theoremmain} are as follows.
Let $m_{\nu_{j}} ( z ) $ and $m_{\nu_{1}\boxplus\nu_{2}}(z)$
denote the Stieltjes transforms of $\nu_{j}$ and $\nu_{1}\boxplus\nu_{2}$,
respectively, and let $\omega_{\nu,j}$ denote the subordination
functions for the pair $ ( \nu_{1},\nu_{2} )$. First, we prove
that the smallness of $d_{L}(\mu_{j},\nu_{j})$ implies that the
differences $\llvert m_{\nu_{j}}-m_{\mu_{j}}\rrvert$ are small,
and that the differences between the derivatives of $m_{\nu_{j}}$ and $
m_{\mu_{j}}$ are also small. Then we show that this fact, together with
system (\ref{systemfort}), implies that the differences between the
corresponding subordination functions are small. At this stage we need the
assumption of smoothness. Finally, we check that if both the Stieltjes
transforms and the subordination functions of pairs $(\mu_{1},\mu_{2})$
and $(\nu_{1},\nu_{2})$ are close to each other, then the Stieltjes
transforms of $\mu_{1}\boxplus\mu_{2}$ and $\nu_{1}\boxplus\nu_{2}$
are close to each other uniformly on the half-line $\Re z=x$, $\Im z>0$.
This fact implies that the densities of $\mu_{1}\boxplus\mu_{2}$ and
$\nu_{1}\boxplus\nu_{2}$ at $x$ are close to each other.

Before discussing applications of Theorem~\ref{theoremmain}, let us mention
some results which are helpful in checking the assumptions of this theorem.

%
\begin{proposition}
\label{propositionsmoothnessstability}If $\mu\boxplus\mu$ is
(Lebesgue) absolutely continuous in a neighborhood of $x$, and the
density of
$\mu\boxplus\mu$ is positive at $x$, then $\Im\omega_{j}(x)>0$
for $j=1,2$.
\end{proposition}

Another important case is when one of the probability measures has the
semicircle distribution with the density $f_{\mathrm{sc}} ( x )
=\frac{1}{%
2\pi}\sqrt{ ( 4-x^{2} )_{+}}$. Since such a measure, $\mu_{\mathrm{sc}}$,
is absolutely continuous, $\mu_{\mathrm{sc}}\boxplus\mu$ is also
absolutely continuous, for an arbitrary $\mu$.

%
\begin{proposition}
\label{propositionsemicirclestability} If the density of $\mu
_{\mathrm{sc}}\boxplus\mu$ is positive at $x$, and
\[
\bigl\llvert m_{\mu
_{\mathrm{sc}}\boxplus
\mu} ( x ) \bigr\rrvert\neq1,
\]
then $\Im\omega_{j}(x)>0$ for $j=1,2$.
\end{proposition}

The proofs of Propositions~\ref{propositionsmoothnessstability} and
\ref{propositionsemicirclestability} will be given in Section
\ref{sectionproofpropostionstability}.

Now let us turn to applications. Theorem~\ref{theoremmain} can be applied
to derive some old and new results about sums of free random variables and
about eigenvalues of large random matrices.\vadjust{\goodbreak}

Recall that if $X_{1},\ldots,X_{n}$ are free, identically distributed
self-adjoint random variables with finite variance $\sigma^{2}$, then
\cite{voiculescu83,maassen92} $S_n:=(X_1+\cdots+X_n)/\break(\sigma\sqrt
{n}) $
converges in distribution to a random variable $X$ with the standard
semicircle law.

In terms of free convolutions, it means that if $\mu$ is a probability
measure with variance $\sigma^{2}$, and if
\[
\mu_{n} ( dx ):=\underset{n\ \mathrm{times}} {\underbrace{\mu
\boxplus\cdots\boxplus\mu}}\bigl(\sigma\sqrt{n}\,dx\bigr),
\]
then $\mu_{n}\rightarrow\mu_{\mathrm{sc}}$.

Bercovici and Voiculescu in~\cite{bercovicivoiculescu95} showed that the
convergence in this limit law holds in a stronger sense. Namely,
assuming in
addition that support of $\mu$ is bounded, they showed that $\mu_{n}$ has
a density for all sufficiently large $n$ and that the sequence of these
densities converges uniformly to the density of the semicircle law.
Recently, this result was generalized in~\cite{wang10} to the case of
$\mu_{n}$ with unbounded support and finite variance. Results in
\cite{bercovicivoiculescu95} and~\cite{wang10} can be considered as local
limit versions of the free CLT.

In the first application (Theorem~\ref{theoremlocalCLT}), we give a short
proof of the easier part of the results in \cite
{bercovicivoiculescu95} and
\cite{wang10} by using Theorem~\ref{theoremmain}. (A more difficult
part of
these results concerns the uniformity of the convergence on $\mathbb{R}$.)

In the second application (Theorem~\ref{theoremlocalPoisson}), we
prove an
analogous local limit result for the sums $S_{n}=X_{1,n}+\cdots+X_{n,n}$,
where $X_{i,n}$ are free projection operators with parameters $p_{i,n}$ such
that $\sum_{i=1}^{n}p_{i,n}\rightarrow\lambda$ and $\max_{i}p_{i,n}%
\rightarrow0$ as $n\rightarrow\infty$. The classical analogue of this
situation is the sum of independent indicator random variables, and the
classical result states that the sums converge in distribution to the
Poisson random variable with parameter $\lambda$. A local version of this
result is absent in the classical case because the Poisson random variable
is discrete, and it does not make sense to talk about convergence of
densities. In the free probability case, the limit of the spectral
distributions of $S_{n}$ is the Marchenko--Pastur distribution, which is
absolutely continuous with bounded density for $\lambda>1$. We show
that in
this case the spectral measures of $S_{n}$ have a density for all
sufficiently large $n$ and that the sequence of these densities converges
uniformly to the density of the Marchenko--Pastur law.

In the third application (Theorem~\ref{theoremlocallimitstable}),
we show
that a similar local limit result holds for sums of free $\boxplus$-stable
random variables.

The fourth application (Theorem~\ref{theoremlocallawRMT}) is of a
different kind and is concerned with eigenvalues of large random matrices.
Let $H_{N}=A_{N}+U_{N}B_{N}U_{N}^{\ast}$, where $A_{N}$ and $B_{N}$
are $N$-by-$N$ Hermitian matrices, and $U_{N}$ is a random unitary
matrix with the
Haar distribution on the unitary group $\mathcal{U} ( N )
$. Let $%
\lambda_{1}^{ ( A ) }\geq\cdots\geq\lambda_{N}^{ (
A ) }$ be the eigenvalues of $A_{N}$. Similarly, let $\lambda_{k}^{ ( B
) }$ and $\lambda_{k}^{ ( H ) }$ be ordered
eigenvalues of matrices $B_{N}$ and $H_{N}$, respectively. Define the
\textit{spectral point measures} of $A_{N}$ by $\mu_{A_{N}}:=N^{-1}\sum
_{k=1}^{N}\delta_{\lambda_{k}^{(A)}(H)}$, and define
the spectral point measures of $B_{N}$ and $H_{N}$ similarly.

Assume that $\mu_{A_{N}}\rightarrow\mu_{\alpha}$ and $\mu
_{B_{N}}\rightarrow\mu_{\beta}$, and that the support of $\mu_{A_{N}}$
and $\mu_{B_{N}}$ is uniformly bounded. Let the pair $(\mu_{\alpha
},\mu_{\beta})$ be smooth at $x$.

Define $\mathcal{N}_{I}:=N\mu_{H_{N}} ( I ) $, the number of
eigenvalues of $H_{N}$ in interval $I$, and let $\mathcal{N}_{\eta
} (
x ):=\mathcal{N}_{ ( x-\eta,x+\eta] }$. Finally, assume
that $\eta=\eta( N ) $ and $\frac{1}{\sqrt{\log
(N)}}\ll\eta
(N)\ll1$.

Then, by using the author's previous results from~\cite{kargin10}, and
Theorem~\ref{theoremmain}, it is shown that
\[
\frac{\mathcal{N}_{\eta}(x)}{\eta N}\rightarrow f_{\mu_{\alpha
}\boxplus
\mu_{\beta}}(x)
\]
with probability 1, where $f_{\mu_{\alpha}\boxplus\mu_{\beta}}$ denotes
the density of $\mu_{\alpha}\boxplus\mu_{\beta}$. This result
generalizes the main result in~\cite{pasturvasilchuk00} where it was proved
that $\mu_{H_{N}}\rightarrow\mu_{\alpha}\boxplus\mu_{\beta}$. It can
be interpreted as a local limit law for eigenvalues of a sum of random
Hermitian matrices.

The rest of the paper is organized as follows. Section
\ref{sectionproofmaintheorem} is concerned with the proof of the main
theorem, Section~\ref{sectionproofpropostionstability} contains
proofs of
Propositions~\ref{propositionsmoothnessstability} and
\ref{propositionsemicirclestability}, Section~\ref{sectionapplications}
contains applications, and Section~\ref{sectionconclusion} concludes.

\section{\texorpdfstring{Proof of Theorem \protect\ref{theoremmain}}{Proof of Theorem 1.3}}
\label{sectionproofmaintheorem}

Let $F_{\mu} ( x ) $ and $F_{\nu} ( x ) $
denote the
cumulative distribution functions of the measures $\mu$ and $\nu$,
respectively.

%
\begin{lemma}
\label{lemmacloseness}Suppose that $d_{L} ( \mu,\nu) =s$.
Assume that $h ( x ) $ is a $C^{1}$ real-valued function,
such that
$\int_{-\infty}^{\infty}\llvert h ( u ) \rrvert
\,du<\infty$ and $%
\int_{-\infty}^{\infty}\llvert h^{\prime} ( u ) \rrvert
\,du<\infty
$. Assume in addition that $h ( u ) $ has a finite number of zeros.
Then,
%
%
\begin{equation}\label{definitionDelta}
\Delta:= \int_{\mathbb{R}}\bigl\llvert h ( u ) \bigl[
F_{\nu
} ( \eta u ) -F_{\mu} ( \eta u ) \bigr] \bigr\rrvert\,du
\leq cs\max\bigl\{ 1,\eta^{-1} \bigr\},
\end{equation}
where $c>0$ depends only on $h$.
\end{lemma}

\begin{pf}
Since $h$ is a continuous function with a finite
number of
zeros, we can decompose the set on which $h(u)$ is nonzero into a finite
number of intervals $I_k$ on which $h(u)$ has a constant sign. Note
that it
suffices to estimate the integral on each of these intervals. Consider the
case when $h(u)>0$ on an interval $I_k$. The treatment of the case $h(u)<0$
is similar.

By using the definition of the L\'{e}vy distance, we obtain the following
estimate:
\begin{eqnarray*}
&&
\bigl\llvert F_{\nu} ( \eta u ) -F_{\mu} ( \eta u ) \bigr
\rrvert\\
&&\qquad\leq\max\bigl\{F_{\mu} ( \eta u+s ) -F_{\mu
} ( \eta u
),F_{\nu} ( \eta u+s ) -F_{\nu} ( \eta u ),
\\
&&\hspace*{60pt}F_{\mu} ( \eta u ) -F_{\mu} ( \eta u-s ),F_{\nu
} (
\eta u ) -F_{\nu} ( \eta u-s ) \bigr\}+s.
\end{eqnarray*}
It suffices to estimate
\[
\int_{I_{k}}h ( u ) \bigl\{F_{\mu} ( \eta u+s )
-F_{\mu
} ( \eta u ) +s\bigr\}\,du,
\]
since the other cases are similar.

First of all, note that
%
%
\begin{equation}\label{estimate0}
\int_{I_{k}}h ( u ) s\,du\leq s\int_{-\infty}^{\infty}
\bigl\llvert h ( u ) \bigr\rrvert\,du\leq cs.
\end{equation}

Next, let $\widetilde{I}_{k}=I_{k}+s/\eta$. Then,
\[
\int_{I_{k}}h ( u ) F_{\mu} ( \eta u+s ) \,du=\int
_{%
\widetilde{I}_{k}}h ( t-s/\eta) F_{\mu} ( \eta t ) \,dt
\]
and therefore,
%
%
\begin{eqnarray}\label{estimate1}
&&
\int_{I_{k}}h ( u ) \bigl[ F_{\mu} ( \eta u+s )
-F_{\mu
} ( \eta u ) \bigr] \,du \nonumber\\
&&\qquad\leq\int_{I_{k}\cap\widetilde
{I}_{k}}%
\bigl[ h ( t-s/\eta) -h ( t ) \bigr] F_{\mu
} ( \eta t ) \,dt
\\
&&\qquad\quad{}+\int_{I_{k}\bigtriangleup\widetilde{I}_{k}}\max\bigl( \bigl\llvert
h ( t-s/\eta) \bigr
\rrvert,\bigl\llvert h ( t ) \bigr\rrvert\bigr) F_{\mu} ( \eta t ) \,dt.
\nonumber
\end{eqnarray}

For the first integral in this estimate, we can use the fact that
\[
h ( t-s/\eta) -h ( t ) =-\int_{t-s/\eta
}^{t}h^{\prime
}
( \xi) \,d\xi
\]
and therefore,
%
%
\begin{eqnarray}\label{estimate2}\quad
\biggl\llvert\int_{I_{k}\cap\widetilde{I}_{k}} \bigl[ h ( t-s/\eta)
-h ( t )
\bigr] F_{\mu} ( \eta t ) \,dt\biggr\rrvert&\leq&\int
_{\mathbb{R}}\int_{t-s/\eta}^{t}\bigl\llvert
h^{\prime} ( \xi) \bigr\rrvert F_{\mu} ( \eta t ) \,d\xi\,dt
\nonumber\\
&=&\int_{\mathbb{R}}\bigl\llvert h^{\prime} ( \xi) \bigr
\rrvert\biggl( \int_{\xi}^{\xi+s/\eta}F_{\mu} (
\eta t ) \,dt \biggr) \,d\xi
\\
&\leq&\frac{s}{\eta}\int_{\mathbb{R}}\bigl\llvert
h^{\prime} ( \xi) \bigr\rrvert\,d\xi.
\nonumber
\end{eqnarray}

For the second integral, we note that
%
%
\begin{eqnarray}\label{estimate3}
\int_{I_{k}\bigtriangleup\widetilde{I}_{k}}\max\bigl( \bigl\llvert h (
t-s/\eta) \bigr
\rrvert,\bigl\llvert h ( t ) \bigr\rrvert\bigr) F_{\mu} ( \eta t ) \,dt &
\leq&\sup\bigl\llvert h ( t ) \bigr\rrvert\llvert I_{k}\bigtriangleup
\widetilde{I}_{k}%
\rrvert
\nonumber\\[-8pt]\\[-8pt]
&\leq&2\sup\bigl\llvert h ( t ) \bigr\rrvert s/\eta.
\nonumber
\end{eqnarray}

By using estimates (\ref{estimate0}), (\ref{estimate1}), (\ref{estimate2})
and (\ref{estimate3}), we obtain
\[
\Delta\leq cs\max\bigl\{ 1,\eta^{-1} \bigr\},
\]
where $c$ depends only on function $h$.\vadjust{\goodbreak}
\end{pf}

Now, let $m_{\mu} ( z ) $ and $m_{\nu} ( z ) $ denote
the Stieltjes transforms of the probability measures $\mu$ and $\nu$,
respectively.

%
\begin{lemma}
\label{lemmaclosenessSttransforms}Let $d_{L} ( \mu,\nu
) =s$
and $z=x+i\eta$, where $\eta>0$. Then:

\begin{longlist}[(a)]
\item[(a)]
$\llvert m_{\mu} ( z ) -m_{\nu} ( z )
\rrvert
<cs\eta^{-1}\max\{ 1,\eta^{-1} \} $ where $c$ is a positive
constant, and

\item[(b)] $\llvert\frac{d^{r}}{dz^{r}} ( m_{\mu} ( z )
-m_{\nu
} ( z ) ) \rrvert<c_{r}s\eta^{-1-r}\max
\{ 1,\eta^{-1} \} $ where $c_{r}$ are positive constants.
\end{longlist}
\end{lemma}

\begin{pf}
(a) By integration by parts,
\[
m_{\mu} ( z ) =\int_{\mathbb{R}}\frac{F_{\mu} (
\lambda
) }{ ( \lambda-z )^{2}}\,d\lambda.
\]
Hence, setting $u= ( \lambda-x ) /\eta$,
\begin{eqnarray*}
\Im m_{\mu} ( z ) &=& \frac{2}{\eta}\int_{\mathbb
{R}}F_{\mu}
( x+\eta u ) \frac{u\,du}{ ( 1+u^{2} )^{2}},
\\
\Re m_{\mu} ( z ) &=& \frac{1}{\eta}\int_{\mathbb
{R}}F_{\mu}
( x+\eta u ) \frac{ ( u^{2}-1 ) \,du}{ (
1+u^{2} )^{2}},
\end{eqnarray*}
and similar formulas hold for $\Im m_{\nu} ( z ) $ and $\Re
m_{\nu
} ( z )$. Since $u ( 1+u^{2} )^{-2}$ and $ (
u^{2}-1 ) ( 1+u^{2} )^{-2}$ satisfy the assumptions
of Lemma %
\ref{lemmacloseness}, Claim (a) follows. Claim (b) can be derived
similarly by writing
\begin{eqnarray*}
\frac{d^{r}}{dz^{r}}m_{\mu} ( z ) &=& ( r+1 ) !\int_{%
\mathbb{R}}
\frac{F_{\mu} ( \lambda) \,d\lambda}{ (
\lambda
-x-i\eta)^{r+2}}
\\
&=&\frac{ ( r+1 ) !}{\eta^{r+1}}\int_{\mathbb{R}}\frac
{1}{ (
u-i )^{r+2}}F_{\mu}
( \eta u+x ) \,du,
\end{eqnarray*}
separating imaginary and real parts of the integrand, and applying
Lemma~\ref{lemmacloseness}.
\end{pf}

%
\begin{lemma}
\label{lemmainitialerror} Assume that the pair $(\mu_{1},\mu_{2})$ is
smooth at $x$. Suppose that $(\nu_{1},\nu_{2})$ is another pair of
probability measures such that $d_{L}(\mu_{j},\nu_{j})<s$ for
$j=1,2$. Let
$\Re z=x$ and $\Im z\geq0$. Then
\[
\biggl\llvert\frac{1}{z-\omega_{\mu,1} ( z ) -\omega_{\mu,2} ( z )
}-m_{\nu_{j}} \bigl( \omega_{\mu,j} (
z ) \bigr) \biggr\rrvert\leq c_{\mu}s
\]
for $j=1,2$. Here $c_{\mu}>0$ depends only on $(\mu_{1},\mu_{2})$
and $x$.
\end{lemma}

That is, if we substitute $\omega_{\mu,j} ( z ) $ in the system
for $\omega_{\nu,j} ( z )$, then the equalities will be
satisfied up to a quantity of order $s$.

\begin{pf*}{Proof of Lemma~\ref{lemmainitialerror}}
The functions $\omega_{\mu,j} ( z ) $ satisfy
equations (\ref{systemfort}), which implies that it is enough to
show that
\[
\bigl\llvert m_{\nu_{j}} \bigl( \omega_{\mu,j} ( z ) \bigr)
-m_{\mu_{j}} \bigl( \omega_{\mu,j} ( z ) \bigr) \bigr\rrvert<cs
\]
for $j=1,2$. Note that $\min_{j=1,2}\{\Im( \omega_{\mu,j} (
x ) ) \}>0$ by the assumption of smoothness of $ (
\mu_{1},\mu_{2} )$. In addition, $\Im( \omega_{\mu,j} (
x+i\eta) ) \geq\eta$ for all $\eta>0$. Hence, by continuity
of $\omega_{\mu,j} ( x+i\eta) $ in $\eta$, we have
$\kappa_{j}:=\inf_{\eta\geq0}\omega_{\mu,j} ( x+i\eta) >0$.
Then, by
Lemma~\ref{lemmaclosenessSttransforms},
\[
\bigl\llvert m_{\nu j} \bigl( \omega_{\mu,j} ( z ) \bigr)
-m_{\mu
j} \bigl( \omega_{\mu,j} ( z ) \bigr) \bigr\rrvert<cs\min
\bigl\{ \kappa_{j}^{-1},\kappa_{j}^{-2}
\bigr\}.
\]
\upqed
\end{pf*}

%
\begin{proposition}
\label{propositiondifferenceoft}Assume that a pair of probability
measures $ ( \mu_{1},\mu_{2} ) $ is smooth at $x$. Then
there are
some $s_{\mu,0}>$ and $c_{\mu}>0$ that depend only on $(\mu_{1},\mu_{2})$
and $x$, such that for all pairs of probability measures $ ( \nu_{1},\nu
_{2} ) $ with $d_{L} ( \mu_{j},\nu_{j} ) <
s\leq
s_{\mu,0}$ for $j=1,2$, the limits $\omega_{\nu,j} ( x )
:=\lim_{\eta\downarrow0}\omega_{\nu,j} ( x+i\eta) $ exist,
and it is true that%
\[
\bigl\llvert\omega_{\nu,j} ( x ) -\omega_{\mu,j} ( x ) \bigr
\rrvert\leq c_{\mu}s
\]
for $j=1,2$.
\end{proposition}

%
\begin{corollary}
Assume that the assumptions of Proposition~\ref{propositiondifferenceoft}
hold and that $d_{L} ( \mu_{j},\nu_{j} ) <s\leq s_{\mu,0}$
for $j=1,2$. Then, $\nu_{1}\boxplus\nu_{2}$ is absolutely continuous
in a
neighborhood of $x$, and
\[
\bigl|f_{\nu_{1}\boxplus\nu_{2}} ( x ) -f_{\mu_{1}\boxplus
\mu
_{2}}(x)\bigr|<c_{\mu}s,
\]
where $f_{\mu_{1}\boxplus\mu_{2}}$ and $f_{\nu_{1}\boxplus\nu_{2}}$
denote the densities of $\mu_{1}\boxplus\mu_{2}$ and $\nu_{1}\boxplus
\nu_{2}$, respectively.
\end{corollary}

\begin{pf}
Since $m_{\nu_{1}\boxplus\nu_{2}} (
z ) = ( z-\omega_{\nu,1} ( z ) -\omega_{\nu,2} (
z ) )^{-1}$, Proposition~\ref{propositiondifferenceoft}
implies that the limit $m_{\nu_{1}\boxplus\nu_{2}} ( x )
:=\lim_{\eta\downarrow0}m_{\nu_{1}\boxplus\nu_{2}} ( x+i\eta
) $ exists and
%
%
\begin{equation}\label{closenessStieltjes}
\bigl\llvert m_{\nu_{1}\boxplus\nu_{2}} ( x ) -m_{\mu
_{1}\boxplus
\mu_{2}} ( x ) \bigr\rrvert
<c_{\mu}s.
\end{equation}
By~\cite{belinschi08}, $\nu_{1}\boxplus\nu_{2}$ has no singular
component. Hence, inequality (\ref{closenessStieltjes}) and the absolute
continuity of $\mu_{1}\boxplus\mu_{2}$ in a neighborhood of $x$ imply
that for all sufficiently small~$s$, the measure $\nu_{1}\boxplus\nu_{2}$
is absolutely continuous in a neighborhood of $x$ with the density
$f_{\nu
_{1}\boxplus\nu_{2}} ( x ) =\pi^{-1}\Im( m_{\nu
_{1}\boxplus\nu_{2}} ( x ) ) $, and
\[
\bigl\llvert f_{\nu_{1}\boxplus\nu_{2}} ( x ) -f_{\mu
_{1}\boxplus
\mu_{2}} ( x ) \bigr\rrvert
<c_{\mu}s.
\]
\upqed
\end{pf}

\begin{pf*}{Proof of Proposition~\ref{propositiondifferenceoft}}
Let $F ( \omega)\dvtx\mathbb{C}^{2}\rightarrow\mathbb{C}^{2}$
be defined
by the formula%
\[
F\dvtx\pmatrix{\omega_{1}
\cr
\omega_{2}}
\rightarrow
\pmatrix{ {c} ( z-\omega_{1}-\omega_{2}
)^{-1}-m_{\nu_{1}} ( \omega_{1} )
\cr
( z-\omega_{1}-\omega_{2} )^{-1}-m_{\nu_{2}}
( \omega_{2} )}.
\]
Let us use the norm $\llVert( x_{1},x_{2} ) \rrVert
= ( \llvert x_{1}\rrvert^{2}+\llvert x_{2}\rrvert^{2} )^{1/2}$. By
Lemma~\ref{lemmainitialerror}, $\llVert
F (
\omega_{\mu,1} ( z )$, $\omega_{\mu,2} ( z )
)
\rrVert\leq c_{\mu}s$ for all $z=x+i\eta$ and $\eta\geq0$.

The derivative of $F$ with respect to $\omega$ is
\[
F^{\prime}= \pmatrix{
( z-
\omega_{1}-\omega_{2} )^{-2} & ( z-
\omega_{1}-\omega_{2} )^{-2}-m_{\nu_{1}}^{\prime}
( \omega_{1} )
\cr
( z-\omega_{1}-\omega_{2} )^{-2}-m_{\nu_{2}}^{\prime
}
( \omega_{2} ) & ( z-\omega_{1}-\omega_{2}
)^{-2}}.
\]
The determinant of this matrix is
\[
\bigl[ m_{\nu_{1}}^{\prime} ( \omega_{1} )
+m_{\nu
_{2}}^{\prime
} ( \omega_{2} ) \bigr] ( z-
\omega_{1}-\omega_{2} )^{-2}-m_{\nu_{1}}^{\prime}
( \omega_{1} ) m_{\nu
_{2}}^{\prime
} (
\omega_{2} ).
\]
By the assumption of smoothness and by Lemma
\ref{lemmaclosenessSttransforms}, this is close (i.e., the difference $<cs$
for some $c>0$) to
\[
\bigl[ m_{\mu_{1}}^{\prime} ( \omega_{1} )
+m_{\mu
_{2}}^{\prime
} ( \omega_{2} ) \bigr] ( z-
\omega_{1}-\omega_{2} )^{-2}-m_{\mu_{1}}^{\prime}
( \omega_{1} ) m_{\mu
_{2}}^{\prime
} (
\omega_{2} )
\]
at $ ( \omega_{1},\omega_{2} ) = ( \omega_{\mu,1} (
z ),\omega_{\mu,2} ( z ) ) $ for all
$z=x+i\eta$
with $\eta\geq0$. The latter expression is nonzero by (\ref{genericity}).
In addition, the assumption of smoothness shows that $ ( z-\omega_{\mu
,1} ( z ) -\omega_{\mu,2} ( z ) )^{-2}$ is
bounded for $z=x+i\eta$ with $\eta\geq0$. Hence, the entries of the
matrix $ [ F^{\prime} ]^{-1}$ are bounded at $ (
\omega_{\mu,1} ( z ),\omega_{\mu,2} ( z ) ) $, and
the bound
does not depend on $\eta$. This shows that the operator norm of $ [
F^{\prime} ]^{-1}$ is bounded uniformly in~$\eta$.

Similarly, an explicit calculation of $F^{\prime\prime}$, the assumption
of smoothness of $ ( \mu_{1},\mu_{2} ) $ and Lemma
\ref{lemmaclosenessSttransforms} imply that for all $z=x+i\eta$ with
$\eta
\geq0$, the operator norm of $F^{\prime\prime}$ is bounded
(uniformly in $%
\eta$) for all $ ( \omega_{1},\omega_{2} ) $ in a neighborhood
of $ ( \omega_{\mu,1} ( z ),\omega_{\mu,2} (
z )
) $.

It follows by the Newton--Kantorovich theorem~\cite{kantorovich48}
that if $%
s=\max_{j}d_{L}(\mu_{j},\break\nu_{j})$ is sufficiently small, then the solution
of the equation $F(\omega)=0$ exists for all $z$ with $\Re z=x$ and
$\Im
z\geq0$.

This solution must be $ ( \omega_{\nu,1} ( z )
,\omega_{\nu,2} ( z ) ) $ by the following argument from~\cite{belinschibercovici07}.
A~solution of equation $F(\omega)=0$
satisfies the
following pair of equations:
\begin{eqnarray*}
\omega_{1} &=&z+h_{2} ( \omega_{2} ),
\\
\omega_{2} &=&z+h_{1} ( \omega_{1} ),
\end{eqnarray*}
where
\[
h_{j} ( \omega) =-\omega-\frac{1}{m_{\nu_{j}} (
\omega
) }.
\]
Note in particular that $\Im h_{j} ( \omega) \geq0$ for
all $%
\omega\in\mathbb{C}^{+}$; see, for example,~\cite{bercovicivoiculescu93} or~\cite{maassen92}.

Hence, $\omega_{1}$ is a fixed point of the function
\[
f_{z} ( \omega) =z+h_{2} \bigl( z+h_{1} (
\omega) \bigr),
\]
which maps $\mathbb{C}^{+}$ to $\mathbb{C}^{+}$. For every $z\in
\mathbb{C}%
^{+}$, the function $f_{z} ( \omega) $ is not a conformal
automorphism because it maps $\mathbb{C}^{+}$ to a subset of $\mathbb
{C}%
^{+}+\Im z$, which is a proper subset of $\mathbb{C}^{+}$. In
addition, it is
analytic as a function of $z$ and $\omega$ that maps $\mathbb
{C}^{+}\times
\mathbb{C}^{+}$ to $\mathbb{C}^{+}$. Hence, by Theorem 2.4 in
\cite{belinschibercovici07}, for every $z\in\mathbb{C}^{+}$ the
function $
f_{z} ( \omega) $ has a unique fixed point $\omega_{1} (
z )$.

A similar argument holds for $\omega_{2} ( z ) $, and we conclude
that equation $F(\omega)=0$ has a unique solution in $\mathbb
{C}^{+}\times
\mathbb{C}^{+}$, which necessarily coincides with $ ( \omega_{\nu,1} (
z )$, $\omega_{\nu,2} ( z ) )$.

In addition, this solution satisfies the inequalities
%
%
\begin{equation}\label{inequalitysubordinationfunctions}
\bigl\llvert\omega_{\nu,j} ( z ) -\omega_{\mu,j} ( z ) \bigr
\rrvert<c_{\mu}s,\qquad j=1,2,
\end{equation}
for all $z$ with $\Re z=E$ and $\Im z>0$.

By Theorem 3.3 in~\cite{belinschi08}, the limits
\[
\omega_{\nu,j} ( E ):=\lim_{\eta\downarrow0}\omega_{\nu,j} ( x+i\eta
)
\]
and
\[
\omega_{\mu,j} ( E ):=\lim_{\eta\downarrow0}\omega_{\mu,j} ( x+i\eta
)
\]
exist, and by taking the limits in (\ref
{inequalitysubordinationfunctions}%
), we find that
\[
\bigl\llvert\omega_{\nu,j} ( x ) -\omega_{\mu,j} ( x ) \bigr
\rrvert\leq cs.
\]
\upqed
\end{pf*}

\section{\texorpdfstring{Proofs of Propositions \protect\ref{propositionsmoothnessstability} and \protect\ref{propositionsemicirclestability}}
{Proofs of Propositions 1.4 and 1.5}}
\label{sectionproofpropostionstability}

Recall that a function $f ( x ) $ is said to be H\"{o}lder
continuous at $x_{0}$ if there exist positive constants $\alpha$, $C$
and $%
\varepsilon$ such that $\llvert x-x_{0}\rrvert<\varepsilon$
implies that $\llvert f(x)-f ( x_{0} ) \rrvert
<C\llvert
x-x_{0}\rrvert^{\alpha}$.

%
\begin{lemma}
\label{lemmarealpartSttransform}Suppose that a probability measure
$\mu$
has a density which is positive and H\"{o}lder continuous at $x$. Let $%
m_{\mu} ( z ) $ be the Stieltjes transform of $\mu$. Then $
\llvert m_{\mu} ( x+i\eta) \rrvert\leq M<\infty
$ for
all $\eta>0$.
\end{lemma}

\begin{pf}
The results of Sokhotskyi, Plemelj and Privalov ensure that
the limit of $m_{\mu} ( x+i\eta) $ exists when $\eta
\downarrow
0 $; see Theorems 14.1b and 14.1c in~\cite{henrici86}. In particular this
implies that $m_{\mu} ( x+i\eta) $ is bounded for sufficiently
small $\eta$. In addition, $\llvert m_{\mu} ( x+i\eta)
\rrvert\leq1/\eta$ so it is bounded for large $\eta$. Since
$m_{\mu
} ( x+i\eta) $ is continuous in the upper half-plane, $%
m_{\mu} ( x+i\eta) $ is bounded for all $\eta$, and the
claim of
the lemma follows.
\end{pf}

\begin{pf*}{Proof of Proposition~\ref{propositionsmoothnessstability}}
Note that for the case $\mu_{1}=\mu_{2}=\mu$,
%
%
\begin{equation}\label{formulafort}
\omega_{1} ( z ) =\omega_{2} ( z ) = \bigl(
z-m_{\mu
\boxplus\mu} ( z )^{-1} \bigr) /2.
\end{equation}
Since by assumption $\mu\boxplus\mu$ is absolutely continuous in a
neighborhood of $x$, and its density $f_{\mu\boxplus\mu}$ is
positive at $%
x$, by the results in~\cite{belinschi08} $f_{\mu\boxplus\mu}$ is
analytic and
therefore uniformly H\"{o}lder continuous in a neighborhood of $x$. By
Sokhotskyi, Plemelj and Privalov's results, the limit $m_{\mu\boxplus
\mu
} ( x ) =\lim_{\eta\downarrow0}m_{\mu\boxplus\mu} (
x+i\eta) $ exists and $\Im m_{\mu\boxplus\mu} ( x
) =\pi
f_{\mu\boxplus\mu} ( x ) >0$. Then it follows from\vadjust{\goodbreak}
(\ref{formulafort}) that the limits $\omega_{j} ( x ) =\lim_{\eta
\downarrow0}\omega_{j} ( x+i\eta) $ exist. Moreover, since
\[
\Im\omega_{j} ( z ) =\frac{1}{2} \biggl( \eta+
\frac{\Im
m_{\mu
\boxplus\mu} ( z ) }{\llvert m_{\mu\boxplus\mu
} (
z ) \rrvert^{2}} \biggr)
\]
and by Lemma~\ref{lemmarealpartSttransform}, $\llvert m_{\mu
\boxplus
\mu} ( z ) \rrvert^{2}$ is bounded uniformly in $\eta$,
hence the fact that $\Im m_{\mu\boxplus\mu} ( x ) =\pi
f_{\mu
\boxplus\mu} ( x ) >0$ implies that $\Im\omega_{j} (
x ) >0$. This completes the proof of the proposition.
\end{pf*}

%
\begin{lemma}
\label{lemmasemicircle}If $\mu_{1}$ has the semicircle distribution, then:

\begin{longlist}[(iii)]
\item[(i)] $\omega_{1} ( z ) =z-\omega_{2} ( z )
+ [
z-\omega_{2} ( z ) ]^{-1}$;

\item[(ii)] $m_{\mu_{\mathrm{sc}}\boxplus\mu} ( z ) =\omega_{2} (
z )
-z$;

\item[(iii)] $\omega_{2}(z)$ satisfies the equation
\[
\omega_{2} ( z ) =z+\int\frac{\mu( dx )
}{x-\omega_{2} ( z ) }.
\]
\end{longlist}
\end{lemma}

\begin{pf}
(i) If $\mu_{1}$ has the semicircle distribution,
then $%
m_{\mu_{1}}^{ ( -1 ) }=- ( z+z^{-1} ) $; hence
the first
equation in system (\ref{systemfort}) implies
\[
\omega_{1}=- \biggl( \frac{1}{z-\omega_{1}-\omega_{2}}+z-\omega_{1}-
\omega_{2} \biggr),
\]
which simplifies to
\[
\omega_{1}=z-\omega_{2}+\frac{1}{z-\omega_{2}}.
\]

\mbox{}\hphantom{i}(ii) By using (i),
\[
m_{\mu_{\mathrm{sc}}\boxplus\mu}=\frac{1}{z-\omega_{1}-\omega_{2}}=- ( z-\omega
_{2} ).
\]

(iii) The second equation in system (\ref{systemfort}) becomes%
\[
- \bigl( z-\omega_{2} ( z ) \bigr) =\int\frac{\mu(
dx ) }{x-\omega_{2} ( z ) }.
\]
\upqed
\end{pf}

\begin{pf*}{Proof of Proposition~\ref{propositionsemicirclestability}}
From (ii) in Lemma~\ref{lemmasemicircle},
\[
\Im\omega_{2}(x)=\Im m_{\mu
_{\mathrm{sc}}\boxplus\mu}(x)=\pi f_{\mu_{\mathrm{sc}}\boxplus\mu}(x)>0.
\]
From (i),
\begin{eqnarray*}
\Im\omega_{1} ( x ) &=&\Im\omega_{2} ( x ) \biggl(
-1+%
\frac{1}{\llvert x-\omega_{2}\rrvert^{2}} \biggr)
\\
&=&\Im\omega_{2} ( x ) \biggl( -1+\frac{1}{\llvert
m_{\mu
_{\mathrm{sc}}\boxplus\mu} ( x ) \rrvert^{2}} \biggr).
\end{eqnarray*}
Since $\Im\omega_{2} ( x ) >0$, if $|m_{\mu_{\mathrm{sc}}\boxplus
\mu} ( x ) |^{2}<1$, then $\Im\omega_{1} ( x )
>0$, and
we are done. Two remaining possibilities are $\llvert m_{\mu
_{\mathrm{sc}}\boxplus\mu} ( x ) \rrvert^{2}=1$ and $\llvert
m_{\mu_{\mathrm{sc}}\boxplus\mu} ( x ) \rrvert^{2}>1$.
However, $%
\llvert m_{\mu_{\mathrm{sc}}\boxplus\mu} ( x ) \rrvert^{2}>1$ is
in fact not possible because\vadjust{\goodbreak} this would imply that $\Im\omega_{1} (
x ) <0$, which is ruled out by a general result of Biane. To sum up,
the assumptions $f_{\mu_{\mathrm{sc}}\boxplus\mu}(x)>0$ and $\llvert
m_{\mu
_{\mathrm{sc}}\boxplus\mu} ( x ) \rrvert^{2}\neq1$ imply
that $\Im
\omega_{j} ( x ) >0$.
\end{pf*}

\section{Applications}\label{sectionapplications}

In the first application we re-prove an easier part of the free local limit
theorem which was first demonstrated in~\cite{bercovicivoiculescu95} for
bounded random variables and later generalized in~\cite{wang10} to the case
of unbounded variables with finite variance. We will show the
convergence of
densities, but we will not investigate whether the convergence is
uniform on $%
\mathbb{R}$.

Let $X_{i}$ be a sequence of self-adjoint identically-distributed free
random variables in the sense of free probability theory. Define $%
S_{n}= ( X_{1}+\cdots+X_{n} ) /\sqrt{n}$, and let $\mu$
and $\mu_{n}$ denote the spectral probability measures of $X_{i}$ and $S_{n}$,
respectively. It is known that
\[
\mu_{n} ( dx ) =\underset{n\ \mathrm{times}} {\underbrace{\mu
\boxplus\cdots\boxplus\mu}} \bigl( \sqrt{n}\,dx \bigr).
\]

%
\begin{theorem}
\label{theoremlocalCLT}Suppose $\mu$ has zero mean and unit
variance. Let
$I_{\varepsilon}= [ -2+\varepsilon,2-\varepsilon] $.
Then for
all sufficiently large $n$, $\mu_{n}$ is (Lebesgue) absolutely continuous
everywhere on $I$, and the density $d\mu_{n}/dx$ uniformly converges
on $%
I_{\varepsilon}$ to the density of the standard semicircle law.
\end{theorem}

Note that the results in~\cite{bercovicivoiculescu95} imply that for every
closed interval $J$ outside of $ [ 2,-2 ]$, the measure
$\mu_{n} ( J ) =0$ for all sufficiently large $n$, provided that
$\mu_{1}$ has bounded support. In addition, the uniform convergence on $
I_{\varepsilon}$ can be strengthened to the uniform convergence on
$\mathbb{%
R}$ as in the proof of Theorem~3.4(iii) in~\cite{wang10}.

\begin{pf*}{Proof of Theorem~\ref{theoremlocalCLT}}
Let $\nu_{1,n}$ be the distribution of $ (
X_{1}+\cdots
+X_{ [ n/2 ] } ) /\sqrt{n}$ and $\nu_{2,n}$ be the
distribution of $ ( X_{ [ n/2 ] +1}+\cdots+X_{n}
) /\sqrt{%
n}$. By using the free CLT (Central limit theorem) from~\cite{maassen92}
(which generalizes the free CLT in~\cite{voiculescu83}), we infer that
both $%
\nu_{1,n}$ and $\nu_{2,n}$ converge weakly to $\widetilde{\mu}_{\mathrm{sc}}$,
where $\widetilde{\mu}_{\mathrm{sc}}$ is the semicircle law with variance
$1/2$. It
is easy to calculate for the pair $ ( \widetilde{\mu
}_{\mathrm{sc}},\widetilde{%
\mu}_{\mathrm{sc}} ) $ that%
\[
\omega_{\widetilde{\mu},1}=\omega_{\widetilde{\mu},2}=\frac
{3z+\sqrt{%
z^{2}-4}}{4}
\]
and therefore $\Im\omega_{\widetilde{\mu},j} ( x ) >0$
on $%
I_{\varepsilon}$. (This also follows from Proposition
\ref{propositionsmoothnessstability}.) A~calculation shows that the
genericity condition (\ref{genericity}) is satisfied for each $x\in
I_{\varepsilon}$, and therefore the density of $\nu_{1,n}\boxplus\nu
_{2,n}$ exists for all sufficiently large~$n$, and converges to the density
of $\widetilde{\mu}_{\mathrm{sc}}\boxplus\widetilde{\mu}_{\mathrm{sc}}$ at each
$x\in
I_{\varepsilon}$. A remark after Theorem~\ref{theoremmain} shows
that the
convergence is in fact uniform on $I_{\varepsilon}$. Since $\nu
_{1,n}\boxplus\nu_{2,n}=\mu_{n}$, this implies that the density of
$\mu_{n}$ converges uniformly on $I_{\varepsilon}$ to the density of the
standard semicircle law.
\end{pf*}

In a similar fashion, it is possible to prove the local limit law for the
convergence to the free Poisson distribution.\vadjust{\goodbreak}

Let $ \{ X_{n,i} \}_{i=1}^{n}$ be freely independent self-adjoint
random variables with the distribution $\mu_{n,i}=p_{n,i}\delta_{1}+ (
1-p_{n,i} ) \delta_{0}$. Let $S_{n}=X_{n,1}+\cdots+X_{n,n}$,
and let $%
\mu_{n}$ denote the spectral probability measure of $S_{n}$. Then
\[
\mu_{n} ( dx ) =\mu_{n,1}\boxplus\cdots\boxplus
\mu_{n,n} ( dx ).
\]
Recall that the \textit{Marchenko--Pastur distribution} with parameter
$\lambda
\geq1$ is a probability measure $\mu_{\mathrm{mp}}$ on $\mathbb{R}$, with the
density
\[
f_{\mathrm{mp}}(x)=\frac{\sqrt{4x- ( 1-\lambda+x )^{2}}}{2\pi x}
\]
supported on the interval $ [ x_{\min},x_{\max} ]:=
[ (
1-\sqrt{\lambda} )^{2}, ( 1+\sqrt{\lambda} )^{2} ]$.
In the free probability literature, this distribution is called the
\textit{free Poisson distribution}.

%
\begin{theorem}
\label{theoremlocalPoisson}Assume that $\sum_{i=1}^{n}p_{n,i}\rightarrow
\lambda>1$ and $\max_{i}p_{n,i}\rightarrow0$ as \mbox{$n\rightarrow\infty$}.
Let $I_{\varepsilon}= [ x_{\min}+\varepsilon,x_{\max
}-\varepsilon%
] $. Then for all sufficiently large $n$, $\mu_{n}$ is (Lebesgue)
absolutely continuous everywhere on $I_{\varepsilon}$, and the density
$%
d\mu_{n}/dx$ uniformly converges on $I_{\varepsilon}$ to the density of
the Marchenko--Pastur law with parameter~$\lambda$.
\end{theorem}

The proof of this theorem is similar to the proof of the previous one. The
first step is the weak convergence of $\mu_{n}$. In the case when $%
p_{n,i}=\lambda/n$ for all~$i$, a proof of weak convergence can be
found on
page 34 in~\cite{voiculescudykemanica92}. The general case is a minor
adaptation of this case, and we omit it. Next, we choose $k_{n}$ so that
\[
\sum_{i=1}^{k_{n}}p_{n,i}\leq
\lambda/2<\sum_{i=1}^{k_{n}+1}p_{n,i}
\]
and define $\nu_{1,n}$ and $\nu_{2,n}$ as the spectral probability
measures of $X_{n,1}+\cdots+X_{n,k_{n}}$ and $X_{n,k_{n}+1}+\cdots
+X_{n,n}$, respectively. It is easy to see that both $\nu_{1,n}$ and
$\nu_{2,n}$ converge weakly to $\widetilde{\mu}_{\mathrm{mp}}$, the Marchenko--Pastur
distribution with parameter $\lambda/2$. By using Proposition
\ref{propositionsmoothnessstability}, we conclude that $\Im\omega_{%
\widetilde{\mu}_{\mathrm{mp}},j} ( x ) >0$ on $I_{\varepsilon}$. Moreover,
a~direct calculation shows that
\[
\omega_{\widetilde{\mu},1} ( z ) =\omega_{\widetilde
{\mu},2} ( z ) =\tfrac{1}{4}
\bigl( z+\lambda-1+\sqrt{ \bigl( z- ( 1+\lambda) \bigr)^{2}-4\lambda}
\bigr)
\]
and
\[
m_{\widetilde{\mu}_{\mathrm{mp}}}^{\prime}=\frac{1-\lambda/2}{2z^{2}}+\frac
{%
-z ( 1+\lambda/2 ) + ( 1-\lambda/2 )^{2}}{2z^{2}\sqrt{%
( z- ( 1+\lambda/2 ) )^{2}-2\lambda}}.
\]
After some calculations the genericity condition (\ref{genericity})
can be
simplified to the following inequality:%
%
%
\begin{eqnarray}\label{genericityequationMP}
f ( x,\lambda):\!&=&x^{3}- \bigl( 5+\tfrac{5}{2}\lambda\bigr)
x^{2}+ \bigl( 7+\tfrac{13}{2}\lambda+2\lambda^{2}
\bigr) x\nonumber\\
&&{}-\bigl(3-5\lambda+%
\tfrac{5}{4}\lambda^{2}+
\tfrac{1}{2}\lambda^{3}\bigr)\\
&\neq&0.\nonumber
\end{eqnarray}

%
\begin{figure}

\includegraphics{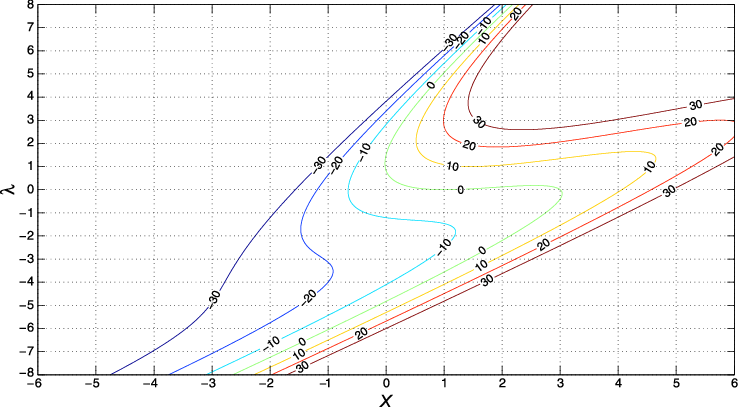}

\caption{Contour plot of the right-hand side of
(\protect\ref{genericityequationMP}).}\label{figgenericity}
\end{figure}

%
\begin{figure}[b]

\includegraphics{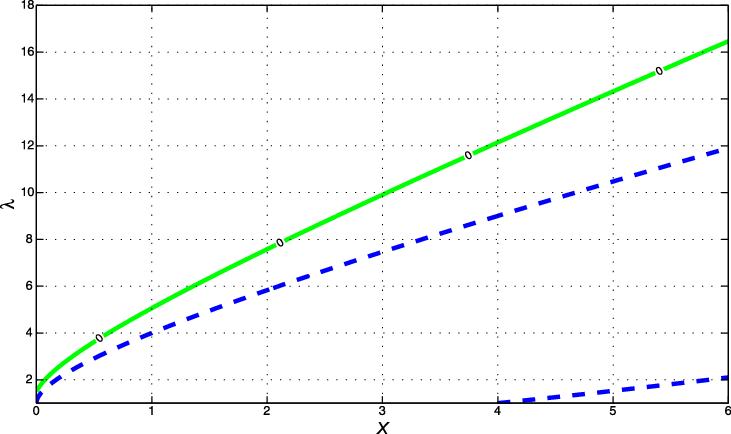}

\caption{The zero set of the right-hand side of
(\protect\ref{genericityequationMP}) compared with the support
bounds for~$x(\lambda)$.}\label{figgenericity2}
\end{figure}

Figure~\ref{figgenericity} shows the contour plot of $f(x,\lambda)$. It
can be seen from this plot and can be checked formally that for
$\lambda>1$, there is only one $x=x ( \lambda) $ that violates
(\ref{genericityequationMP}). Figure~\ref{figgenericity2} shows the
zero set of $f ( x,\lambda) $ for $\lambda>1$, compared with the bounds
on the support of the Marchenko--Pastur distribution. It can be seen
from this graph and can be checked formally that $x(\lambda)<t_{\min} (
\lambda) = ( 1-\sqrt{\lambda} )^{2}$. Consequently, if $x$ is in the
support of
$\widetilde{\mu}_{\mathrm{mp}}\boxplus\widetilde{\mu}_{\mathrm{mp}}$,
the genericity condition (\ref{genericity}) holds, and the pair $ (
\widetilde{\mu}_{\mathrm{mp}},\widetilde{\mu}_{\mathrm{mp}} ) $ is
smooth at $x$. Hence, Theorem~\ref{theoremmain} applies, and the
density of $\mu _{n}=\nu_{1,n}\boxplus\nu_{2,n}$ converges uniformly on
$I_{\varepsilon }$ to the density of
$\widetilde{\mu}_{\mathrm{mp}}\boxplus\widetilde{\mu}_{\mathrm{mp}}$,
that is, to the density of the Marchenko--Pastur distribution with
parameter $\lambda$.

Similar results can be established for other limit theorems, except
that it
is more difficult to check the genericity condition (\ref{genericity})
for a
point in the support of the limit distribution. Here is one more
theorem of
this type. Let measures $\mu$ and $\nu$ be called equivalent ($\mu
\sim
\nu$) if there exist such real $a$ and $b$, with $b>0$, that for every
Borel set $S\subset\mathbb{R}$, $\mu(S)=\nu(bS+a)$. Recall that a measure
$\mu$ is called $\boxplus$-stable, if $\mu\boxplus\mu\sim\mu$. The
measure $\nu$ belongs to the domain of attraction of a $\boxplus$-stable
law~$\mu$, if there exist measures $\nu_{n}$ equivalent to $\nu$ such
that
\[
\underset{n\ \mathrm{times}} {\underbrace{\nu_{n}\boxplus
\nu_{n}\boxplus\cdots\boxplus\nu_{n}}}\rightarrow\mu.
\]
Clearly, in this case there exists a sequence of real constants $b_{n}>0$
and $a_{n}$ such that%
%
%
\begin{equation}\label{coefficientsstablelimittheorem}
\mu_{n}:=\underset{n\ \mathrm{times}} {\underbrace{\nu\boxplus\nu
\boxplus\cdots\boxplus\nu}} ( b_{n}\cdot+a_{n} )
\rightarrow\mu.
\end{equation}
(More about the $\boxplus$-stability of probability measures and its
relation to the classical stability of probability measures can be
found in
\cite{bercovicipata99}.)

%
\begin{theorem}
\label{theoremlocallimitstable}Suppose that a $\boxplus$-stable
distribution $\mu$ is not equivalent to $\delta_{0}$ and that $\nu$
belongs to the domain of attraction of $\mu$. Let $a_{n}$, $b_{n}$ and
$%
\mu_{n}$ be defined as in (\ref{coefficientsstablelimittheorem}),
and let
$J$ be a bounded closed interval such that the density of $\mu$ is strictly
positive on $J$. Then $\mu_{n}$ is (Lebesgue) absolutely continuous on $J$
for all sufficiently large $n$, and there exist such real $\kappa
_{n}>0$ and
$\xi_{n}$ that the density of $\mu_{n} ( \kappa_{n}\cdot+\xi_{n} ) $
converges to the density of $\mu$ at (Lebesgue) almost
all $%
E\in J$.
\end{theorem}

\begin{pf}
Let $J\subset I$, where the inclusion is strict, and
$I$ is a
bounded, closed interval such that density of $\mu$ is strictly
positive on $%
I$. (Interval $I$ exists because by the results of Biane in
\cite{bercovicipata99} $\mu$ is absolutely continuous with analytical density.)

First, note that $\mu_{n}$ is (Lebesgue) absolutely continuous on
$\mathbb{R%
}$ for all sufficiently large $n$. Indeed, for even $n=2k$, the definition
of $\mu_{n}$ implies that $\mu_{2k}=\mu_{k}\boxplus\mu_{k} (
s_{k}^{-1}\cdot-t_{k} ) $ for some constants $t_{k}$ and $s_{k}>0$.
For large $k$, $\mu_{k}$ is close in the L\'{e}vy metric to $\mu$, which
is known to be absolutely continuous. Hence, $\mu_{k}$ has no atoms with
weight $\geq1/2$. This implies that $\mu_{2k}$ has no atoms at all. In
addition, by results of Belinschi, $\mu_{2k}$ has no singular component.
Therefore, $\mu_{2k}$ is absolutely continuous on $\mathbb{R}$ if $k$ is
sufficiently large. The argument for the odd $n=2k+1$ is similar if we write
$\mu_{2k+1}=\mu_{k+1}\boxplus\mu_{k} ( s_{k}\cdot+t_{k}
)$.

In the second and final step, we note that there exists a sequence of
constants $\kappa_{n}>0$ and $\xi_{n}$ such that the density of $\mu
_{n} ( \kappa_{n}\cdot+\xi_{n} ) $ converges to the
density of $%
\mu$ at (Lebesgue) almost all $x\in I$. Indeed, by the stability of
$\mu$, $%
\mu\boxplus\mu=\mu( s\cdot+t ) $ and $\mu$ has positive
analytic density on $I$; therefore, by Proposition
\ref{propositionsmoothnessstability} $\Im\omega_{\mu,j} (
x ) >0$
at all $x\in(I-t)/s$. For almost all points $x$, the genericity
condition (%
\ref{genericity}) holds, since otherwise $k_{\mu} ( x ) $
(in the
genericity condition) would be exactly $0$ which is not possible. For
even $%
n=2k$, we have $\mu_{k}\boxplus\mu_{k}=\mu_{2k} ( s_{k}\cdot
+t_{k} )$, where $s_{k}>0$ and $t_{k}$ are certain real constants.
Hence, by Theorem~\ref{theoremmain} the weak convergence $\mu
_{k}\rightarrow\mu$ implies that the density of $\mu_{k}\boxplus\mu
_{k}\equiv\mu_{2k} ( s_{k}\cdot+t_{k} ) $ converges to the
density of $\mu\boxplus\mu\equiv\mu( s\cdot+t ) $ at almost
all points of $(I-t)/s$. It follows that for $\kappa_{2k}=s/s_{k}>0$
and $%
\xi_{2k}=t-(s/s_{k})t_{k}$, the density of $\mu_{2k} ( \kappa_{2k}\cdot
+\xi_{2k} ) $ converges to the density of $\mu$ at almost
all points of $I$. The case of $\mu_{2k+1}$ can be handled similarly by
considering $\mu_{k}\boxplus\mu_{k+1}$.
\end{pf}

Our next application is of a different kind and answers a question that
arises in the theory of large random matrices.

Let $H_{N}=A_{N}+U_{N}B_{N}U_{N}^{\ast}$, where $A_{N}$ and $B_{N}$
are $N$%
-by-$N$ Hermitian matrices, and $U_{N}$ is a random unitary matrix with the
Haar distribution on the unitary group $\mathcal{U} ( N )$.

Let $\lambda_{1}^{ ( A ) }\geq\cdots\geq\lambda_{N}^{ (
A ) }$ be the eigenvalues of $A_{N}$. Similarly, let $\lambda_{k}^{ ( B
) }$ and $\lambda_{k}^{ ( H ) }$ be ordered
eigenvalues of matrices $B_{N}$ and $H_{N}$, respectively.

Define the\vspace*{1pt} \textit{spectral point measures} of $A_{N}$ by $\mu
_{A_{N}}:=N^{-1}\sum_{k=1}^{N}\delta_{\lambda_{k}^{(A)}(H)}$, and define
the spectral point measures of $B_{N}$ and $H_{N}$ similarly. Let $%
\mathcal{N}_{I}:=N\mu_{H_{N}} ( I ) $ denote the number of
eigenvalues of $H_{N}$ in interval $I$, and let $\mathcal{N}_{\eta
} (
x ):=\mathcal{N}_{ ( x-\eta,x+\eta] }$.

Let the notation $g_{1} ( N ) \ll g_{2}(N)$ mean that $%
\lim_{N\rightarrow\infty}g_{2} ( N ) /g_{1}(N)=+\infty$.

%
\begin{theorem}
\label{theoremlocallawRMT}Assume that:

\begin{longlist}[(4)]
\item[(1)]$\mu_{A_N} \rightarrow\mu_\alpha$ and $\mu_{B_N} \rightarrow
\mu_\beta$;

\item[(2)]$\operatorname{supp}(\mu_{A_N}) \cup\operatorname{supp}(\mu
_{B_N})\subseteq
[-K,K] $ for all $N$;

\item[(3)] the pair $(\mu_{\alpha},\mu_{\beta})$ is smooth at $x$;

\item[(4)]$\frac{1}{\sqrt{\log(N)}} \ll\eta(N) \ll1$.
\end{longlist}

Then
\[
\frac{\mathcal{N}_{\eta}(x)}{2\eta N}\rightarrow f_{\mu_{\alpha
}\boxplus
\mu_{\beta}}(x)
\]
with probability 1, where $f_{\mu_{\alpha}\boxplus\mu_{\beta}}$ denotes
the density of $\mu_{\alpha}\boxplus\mu_{\beta}$.
\end{theorem}

Previously, it was shown by Pastur and Vasilchuk in~\cite{pasturvasilchuk00}
that assumption (1) together with a weaker version of assumption (2) implies
that $\mu_{H_{N}}\rightarrow\mu_{\alpha}\boxplus\mu_{\beta}$ with
probability 1. Theorem~\ref{theoremlocallawRMT} says that the convergence
of $\mu_{H_{N}}$ to $\mu_{\alpha}\boxplus\mu_{\beta}$ holds on the
level of densities, so it can be seen as a local limit law for the
eigenvalues of the sum of random Hermitian matrices.

\begin{pf*}{Proof of Theorem~\ref{theoremlocallawRMT}}
In Theorem 2 in~\cite{kargin10}, it was shown that the
following claim holds. Suppose that $\eta=\eta( N ) $ and
$1/%
\sqrt{\log N}\ll\eta( N ) \ll1$. Assume that the
measure $\mu_{A_{N}}\boxplus\mu_{B_{N}}$ is absolutely continuous, and its
density is
bounded by a constant $T_{N}$. Then, for all sufficiently large $N$,
%
%
\begin{equation}\label{inequalityRMTlocallaw}
P \biggl\{ \sup_{x}\biggl\llvert\frac{\mathcal{N}_{\eta} (
x ) }{%
2N\eta}-f_{\boxplus,N}
( x ) \biggr\rrvert\geq\delta\biggr\} \leq\exp\biggl( -c\delta^{2}
\frac{ ( \eta N )^{2}}{ ( \log
N )^{2}} \biggr),
\end{equation}
where $c>0$ depends only on $K_{N}:=\max\{ \llVert
A_{N}\rrVert,\llVert B_{N}\rrVert\} $ and $T_{N}$. Here
$f_{\boxplus,N}$
denotes the density of $\mu_{A_{N}}\boxplus\mu_{B_{N}}$.

This statement can be modified so that the supremum in the inequality holds
for $x$ in an interval, provided that the density of $\mu
_{A_{N}}\boxplus
\mu_{B_{N}}$ is bounded by a constant $T_{N}$ in the following
interval:%
%
%
\begin{equation}\label{inequalityRMTlocallaw2}
P \biggl\{ \sup_{x\in( a,b ) }\biggl\llvert\frac{\mathcal
{N}_{\eta
} ( x ) }{2N\eta}-f_{\boxplus,N}
( x ) \biggr\rrvert\geq\delta\biggr\} \leq\exp\biggl( -c\delta^{2}
\frac{ ( \eta
N )^{2}%
}{(\log N)^{2}} \biggr).
\end{equation}

Since assumptions (1) and (3) hold, we can use Theorem
\ref{theoremmain} and infer that $f_{\boxplus,N} ( x )
\rightarrow
f_{\mu_{\alpha}\boxplus\mu_{\beta}} ( x )$. In particular,
the sequence of densities $f_{\boxplus,N} ( x ) $ is uniformly
bounded by a constant $T$. This fact and assumption (2) imply that the
positive constant $c$ in (\ref{inequalityRMTlocallaw}) can be chosen
independently of $N$. By using the Borel--Cantelli lemma, we can conclude
that%
\[
\frac{\mathcal{N}_{\eta} ( x ) }{2N\eta}\rightarrow f_{\mu
_{\alpha}\boxplus\mu_{\beta}} ( x )
\]
with probability 1.
\end{pf*}

\section{Conclusion}

\label{sectionconclusion}

We have proved that if probability measures $\nu_{1}$ and $\nu_{2}$ are
sufficiently close to probability measures $\mu_{1}$ and $\mu_{2}$ in
the L%
\'{e}vy distance, and if $\mu_{1}\boxplus\mu_{2}$ is sufficiently smooth
at $x$, then $\nu_{1}\boxplus\nu_{2}$ is absolutely continuous at
$x$, and
its density is close to the density of $\mu_{1}\boxplus\mu_{2}$.

We have applied this result to derive several local limit law results
for sums of
free random variables and for eigenvalues of a sum of random Hermitian
matrices.

\section*{Acknowledgment}
I would like to thank Diana Bloom for her editorial help and an
anonymous referee for helpful suggestions.


%

\printaddresses

\end{document}